\documentclass[12pt]{amsart}
\usepackage{amssymb}

\numberwithin{equation}{section}
\theoremstyle{plain}
\newtheorem{theorem}{Theorem}[section]
\newtheorem{proposition}[theorem]{Proposition}
\newtheorem{corollary}[theorem]{Corollary}
\newtheorem{lemma}[theorem]{Lemma}

\theoremstyle{definition}

\theoremstyle{remark}

\def\R{{\mathbb R}}
\def\C{{\mathbb C}}
\def\Rn{{{\mathbb R}^n}}

\def\Sph{{\mathbb S}}
\def\FT{{\mathcal F}}
\def\L2tx{{L^2(\R_t\times\R^n_x)}}

\def\p#1{{\left({#1}\right)}}
\def\b#1{{\left\{{#1}\right\}}}
\def\br#1{{\left[{#1}\right]}}
\def\n#1{{\left\|{#1}\right\|}}
\def\abs#1{{\left|{#1}\right|}}
\def\jp#1{{\left\langle{#1}\right\rangle}}

\def\supp{\operatorname{supp}}

\def\ka{\kappa}

\title{Smoothing properties of inhomogeneous
equations via canonical transforms}

\author[]{Michael Ruzhansky and Mitsuru Sugimoto}
\address{
  Michael Ruzhansky:
  \endgraf
  Department of Mathematics
  \endgraf
  Imperial College London
  \endgraf
  180 Queen's Gate, London SW7 2AZ, UK
  \endgraf
  {\it E-mail address} {\rm m.ruzhansky@imperial.ac.uk}
  \endgraf
  \medskip
  Mitsuru Sugimoto:
  \endgraf
  Graduate School of Mathematics
  \endgraf
  Nagoya University
  \endgraf
  Furocho, Chikusa-ku, Nagoya 464-8602, Japan
  \endgraf
  {\it E-mail address} {\rm sugimoto@math.nagoya-u.ac.jp}
  }

\thanks{The authors were supported by the Daiwa Anglo-Japanese Foundation.
The first author was also supported by the EPSRC Leadership Fellowship
EP/G007233/1. }

\date{\today}

\begin{document}

\maketitle

\begin{abstract}
The paper describes a new approach to global smoothing problems
for inhomogeneous dispersive evolution equations 
based on an idea of canonical transformation.
In our previous papers \cite{RS1, RS6},
we introduced such a method to show
global smoothing estimates for homogeneous dispersive equations.
It is remarkable that this method allows us to carry out a global
microlocal reduction of equations to some low dimensional model cases.
The purpose of this paper is to pursue 
the same treatment for inhomogeneous equations.
Especially, time-global smoothing estimates for the operator
$a(D_x)$ with lower order terms are the benefit
of our new method.
\end{abstract}

\section{Introduction}
This article consists partly of a survey of the arguments developed
in author's recent paper \cite{RS6} (Sections \ref{SECTION:canonical}
and \ref{SECTION:main}) and partly of obtaining new results via the
extension and continuation of these arguments
(Sections \ref{SECTION:model-inh} and \ref{SECTION:inhomogeneous}).

Let us first consider the following Schr\"odinger equation:
\[
\left\{
\begin{aligned}
\p{i\partial_t+\Delta_x}\,u(t,x)&=0\quad\text{in $\R_t\times\R^n_x$},\\
u(0,x)&=\varphi(x)\quad\text{in $\R^n_x$}.
\end{aligned}
\right.
\]
We know that the solution
operator $e^{it\Delta_x}$ preserves the $L^2$-norm 
for each fixed $t\in\R$.
On the other hand, the extra gain of regularity of order 
$1/2$ in $x$ can
be observed if we integrate the solution in $t$.
For example we have the estimate
\[
\n{\langle x\rangle^{-s}|D_x|^{1/2} e^{it\Delta_x}\varphi}
_{L^2\p{\R_t\times\R^n_x}}
\le C\n{\varphi}_{L^2\p{\R^n_x}} \qquad (s>1/2)
\]
for $u=e^{it\Delta_x}\varphi$
and this estimate was first given by Ben-Artzi and Klainerman 
\cite{BK} ($n\geq3$).
Since the independent pioneering works
by Constantin and Saut \cite{CS},
Sj\"olin \cite{Sj} and Vega \cite{V}, 
the local, then the global smoothing estimates for 
Schor\"odinger or more general dispersive equations have been 
intensively investigated.
(Smoothing for generalised Korteweg-de Vries equations
was already studied by Kato \cite{Ka2}.)
There has already been a lot of
literature on this subject:
Ben-Artzi and Devinatz \cite{BD1, BD2}, Chihara \cite{Ch},
Hoshiro \cite{Ho1, Ho2},
Kato and Yajima \cite{KY},
Kenig, Ponce and Vega \cite{KPV1}--\cite{KPV5}, 
Linares and Ponce \cite{LP}, Simon \cite{Si}, 
Sugimoto \cite{Su1, Su2},
Walther \cite{Wa1, Wa2}, and many others.

In our previous papers \cite{RS1, RS6},
we introduced a new method to show
global smoothing estimates for Schor\"odinger equations,
or more generally, those for homogeneous dispersive equations:
\begin{equation}\label{EQ:Disp}
\left\{
\begin{aligned}
\p{i\partial_t+a(D_x)}\,u(t,x)&=0\quad\text{in $\R_t\times\R^n_x$},\\
u(0,x)&=\varphi(x)\quad\text{in $\R^n_x$}.
\end{aligned}
\right.
\end{equation}
where $a(\xi)$ is a real-valued
function of $\xi=(\xi_1,\ldots,\xi_n)$
with the growth of order $m$,
and $a(D_x)$ is the corresponding operator.
The main idea was to change the equation
\[
\p{i\partial_t+a(D_x)}\,u(t,x)=0
\quad {\rm to}\quad
\p{i\partial_t+\sigma(D_x)}\,v(t,x)=0,
\]
where the operators $a(D_x)$ and $\sigma(D_x)$ are related with each other
by the relation
$$
a(\xi)=\p{\sigma\circ\psi}(\xi).
$$ 
Such an idea can be realised by a canonical transformation
$T$ in the following way:
\[
a(D_x)\circ T=T\circ\sigma(D_x).
\] 
If now operators $T$ and $T^{-1}$ are bounded in $L^2(\R^n_x)$
and in weighted $L^2(\R^n_x)$ respectively,
we can reduce global smoothing estimates for $u=e^{ita(D_x)}\varphi$
to those for $v=e^{it\sigma(D_x)}\varphi$.
It is remarkable that 
the method of canonical transformations described above
allows us to carry out a global microlocal reduction
of equation \eqref{EQ:Disp} to the
model cases $a(\xi)=|\xi_n|^m$ (elliptic case)
or $a(\xi)=\xi_1|\xi_n|^{m-1}$ (non-elliptic case)
under a dispersiveness condition.
\par
The purpose of this paper is to pursue 
the same treatment for inhomogeneous equations:
\[
\left\{
\begin{aligned}
\p{i\partial_t+a(D_x)}\,u(t,x)&=f(t,x)\quad\text{in $\R_t\times\R^n_x$},\\
u(0,x)&=0\quad\text{in $\R^n_x$}.
\end{aligned}
\right.
\]
We will obtain the 
corresponding results on the global smoothing for
solutions to inhomogeneous problems.
There are considerably
less results on this topic available in the literature.
Mostly the Schr\"odinger equation was treated
(e.g. Linares and Ponce \cite{LP}, Kenig, Ponce and
Vega \cite{KPV5}), or the one dimensional case
(Kenig, Ponce and Vega \cite{KPV3, KPV4} or Laurey \cite{La}).
Some more general results on the local
smoothing for dispersive operators were obtained by
Chihara \cite{Ch} and Hoshiro \cite{Ho2}, and for 
dispersive differential operators by Koch and Saut \cite{KoSa}.
In this paper we will extend these results in two directions:
we will establish the global smoothing for rather general
dispersive equations of different orders in all dimensions.
Especially, these kinds of time-global estimate 
for the operator $a(D_x)$ with lower order terms are the benefit
of our new method. The treatment of the inhomogeneous equations
may allow one to treat nonlinear equations with lower order terms and with
corresponding nonlinearities, see the author's paper \cite{RScmp}
for one example.
\par
We will explain the organisation of this paper.
In Section \ref{SECTION:canonical}, we introduce
our main tools established by the authors in \cite{RS6}, which
originate in the idea of canonical transformation.
In Section \ref{SECTION:main}, we list results of smoothing estimates
for homogeneous equations
which were partially announced in \cite{RS1}
and will be completely given in \cite{RS6}.
We also explain how general cases can be reduced
to the model estimates by using canonical transformation.
Sections \ref{SECTION:model-inh} and \ref{SECTION:inhomogeneous} are
devoted to non-homogeneous problems as a counterpart of Section
\ref{SECTION:main}. Model estimates will be given in
Section \ref{SECTION:model-inh}, and estimates for general cases will
be given in Section \ref{SECTION:inhomogeneous}
by using the idea of canonical transformation.
Such argument and related results were partly announced in \cite{RS4}.

\par
Finally we comment on the notation used in this paper.
As usual, we will denote $D_{x_j}=-i\partial_{x_j}$ and
view operators $a(D_x)$ as Fourier multipliers.
Constants denoted by letter $C$ 
in estimates are always positive and
may differ on different 
occasions, but will still be denoted by the same letter.

\section{Canonical transforms}
\label{SECTION:canonical}
We will first review our main tool to reduce general operators to
normal forms discussed in \cite{RS6}.

Let $\psi:\Gamma\to\widetilde{\Gamma}$
be a $C^\infty$-diffeomorphism between open sets
$\Gamma\subset\R^n$ and $\widetilde{\Gamma}\subset\R^n$.
We always assume that
\begin{equation}\label{infty}
C^{-1}\leq\abs{\det \partial\psi(\xi)}\leq C\quad(\xi\in\Gamma),
\end{equation}
for some $C>0$.
We set formally
\begin{align*}
&I_\psi u(x)
=\FT^{-1}\left[\FT u\p{\psi(\xi)}\right](x)
=(2\pi)^{-n}\int_{\R^n}\int_{\R^n}
   e^{i(x\cdot\xi-y\cdot\psi(\xi))}u(y)\, dyd\xi,
\\
&I^{-1}_\psi u(x)
=\FT^{-1}\left[\FT u\p{\psi^{-1}(\xi)}\right](x)
=(2\pi)^{-n}\int_{\R^n}\int_{\R^n}
  e^{i(x\cdot\xi-y\cdot\psi^{-1}(\xi))}u(y)\, dyd\xi.
\end{align*}
The operators $I_\psi$ and $I^{-1}_\psi$ can be justified
by using cut-off functions
$\gamma\in C^\infty(\Gamma)$ and
$\widetilde{\gamma}=\gamma\circ\psi^{-1}\in C^\infty(\widetilde{\Gamma})$
which satisfy $\supp\gamma\subset\Gamma$,
$\supp\widetilde{\gamma}\subset\widetilde{\Gamma}$.
We set
\begin{equation}\label{DefI0}
\begin{aligned}
I_{\psi,\gamma} u(x)
&=\FT^{-1}\left[\gamma(\xi)\FT u\p{\psi(\xi)}\right](x)
\\
&=(2\pi)^{-n}\int_{\R^n}\int_{\Gamma}
 e^{i(x\cdot\xi-y\cdot\psi(\xi))}\gamma(\xi)u(y) dyd\xi,
\\
I_{\psi,\gamma}^{-1} u(x)
&=\FT^{-1}\left[\widetilde{\gamma}(\xi)\FT
u\p{\psi^{-1}(\xi)}\right](x)
\\
&=(2\pi)^{-n}\int_{\R^n}
\int_{\widetilde{\Gamma}}
 e^{i(x\cdot\xi-y\cdot\psi^{-1}(\xi))}
 \widetilde{\gamma}(\xi)u(y) dyd\xi.
\end{aligned}
\end{equation}
In the case that $\Gamma$, $\widetilde{\Gamma}\subset\R^n\setminus0$
are open cones,
we may consider the homogeneous $\psi$ and $\gamma$ which satisfy
$\supp\gamma\cap \Sph^{n-1}\subset\Gamma\cap \Sph^{n-1}$ and
$\supp\widetilde{\gamma}\cap 
\Sph^{n-1}\subset\widetilde{\Gamma}\cap \Sph^{n-1}$,
where $\Sph^{n-1}=\b{\xi\in\Rn: |\xi|=1}$.
Then we have the expressions for compositions
\begin{equation}\label{eq:cut}
I_{\psi,\gamma}=\gamma(D_x)\cdot I_\psi=I_\psi\cdot\widetilde{\gamma}(D_x),\quad
I_{\psi,\gamma}^{-1}
=\widetilde{\gamma}(D_x)\cdot I^{-1}_\psi=I^{-1}_\psi\cdot\gamma(D_x),
\end{equation}
and the identities
\begin{equation}\label{eq:id}
I_{\psi,\gamma}\cdot I_{\psi,\gamma}^{-1}=\gamma(D_x)^2,\quad
I_{\psi,\gamma}^{-1}\cdot I_{\psi,\gamma}=\widetilde{\gamma}(D_x)^2.
\end{equation}
We have also the formulae
\begin{equation}\label{eq:cnon}
I_{\psi,\gamma}\cdot\sigma(D_x)
=\p{\sigma\circ\psi}(D_x)\cdot I_{\psi,\gamma},\quad
I_{\psi,\gamma}^{-1}\cdot \p{\sigma\circ\psi}(D_x)
=\sigma(D_x)\cdot I_{\psi,\gamma}^{-1}.
\end{equation}
\par
\medskip
We also introduce the weighted $L^2$-spaces.
For the weight function $w(x)$, let $L^2_{w}(\R^n;w)$ be
the set of measurable functions $f:\Rn\to\C$ 
such that the norm
\[
\n{f}_{L^2(\R^n;w)}
=\p{\int_{\R^n}\abs{w(x) f(x)}^2\,dx}^{1/2}
\]
is finite.
Then, from the relations \eqref{eq:cut}, \eqref{eq:id}, and \eqref{eq:cnon},
we obtain the following fundamental theorem (\cite[Theorem 4.1]{RS6}):
\par
\medskip
\begin{theorem}\label{Th:reduction}
Assume that the operator $I_{\psi,\gamma}$ defined by \eqref{DefI0}
is $L^2(\R^n;w)$--bounded.
Suppose that we have the estimate
\begin{equation}\label{red}
\n{w(x)\rho(D_x)e^{it\sigma(D_x)}\varphi(x)}_{L^2\p{\R_t\times\R^n_x}}
\leq C\n{\varphi}_{L^2\p{\R^n_x}}
\end{equation}
for all $\varphi$ such
that $\supp\widehat{\varphi}\subset\supp\widetilde\gamma$, where
$\widetilde\gamma=\gamma\circ\psi^{-1}$.
Assume also that the function
\begin{equation}\label{bdd}
q(\xi)=\frac{\gamma\cdot\zeta}{\rho\circ \psi}(\xi)
\end{equation}
is bounded.
Then we have
\begin{equation}\label{org}
\n{w(x)\zeta(D_x)e^{ita(D_x)}\varphi(x)}_{L^2\p{\R_t\times\R^n_x}}
\leq C\n{\varphi}_{L^2\p{\R^n_x}}
\end{equation}
for all $\varphi$
such that $\supp\widehat{\varphi}\subset\supp\gamma$,
where $a(\xi)=(\sigma\circ\psi)(\xi)$.
\end{theorem}
\medskip
\par
Note that $e^{ita(D_x)}\varphi(x)$ and $e^{it\sigma(D_x)}\varphi(x)$ are
solutions to
\[
\left\{
\begin{aligned}
\p{i\partial_t+a(D_x)}\,u(t,x)&=0,\\
u(0,x)&=\varphi(x),
\end{aligned}
\right.
\quad {\rm and}\quad
\left\{
\begin{aligned}
\p{i\partial_t+\sigma(D_x)}\,v(t,x)&=0,\\
v(0,x)&=g(x),
\end{aligned}
\right.
\]
respectively.
Theorem \ref{Th:reduction} means that smoothing estimates for equation with
$\sigma(D_x)$ implies those with $a(D_x)$ if the canonical transformations
which
relate them are bounded on weighted $L^2$-spaces.
The same thing is true for inhomogeneous equations
\[
\left\{
\begin{aligned}
\p{i\partial_t+a(D_x)}\,u(t,x)&=f(t,x),\\
u(0,x)&=0,
\end{aligned}
\right.
\quad {\rm and}\quad
\left\{
\begin{aligned}
\p{i\partial_t+\sigma(D_x)}\,v(t,x)&=f(t,x),\\
v(0,x)&=0,
\end{aligned}
\right.
\]
whose solutions are 
$-i\int^t_0e^{i(t-\tau)a(D_x)}f(\tau,x)\,d\tau$
and
$-i\int^t_0e^{i(t-\sigma)a(D_x)}f(\tau,x)\,d\tau$,
respectively.
The only difference is that we need the weighted $L^2$--boundedness
of the operator $I_{\psi,q}^{-1}$
instead of just the $L^2$--boundedness of it induced by the boundedness
of $q(\xi)$:
\par
\medskip
\begin{theorem}\label{Th:redhom}
Assume that the operator $I_{\psi,\gamma}$ defined by \eqref{DefI0}
is $L^2(\R^n;w)$--bounded.
Suppose that we have the estimate
\[
\n{w(x)\rho(D_x)\int^t_0
e^{i(t-\tau)\sigma(D_x)}f(\tau,x)\,d\tau}_{L^2\p{\R_t\times\R^n_x}}
\leq C\n{v(x)f(t,x)}_{L^2\p{\R_t\times\R^n_x}}
\]
for all $f$ such that $\supp \FT_x f(t,\cdot)\subset\supp\widetilde\gamma$,
where $\widetilde\gamma=\gamma\circ\psi^{-1}$.
Also assume that
the operator
$I_{\psi,q}^{-1}$
defined by \eqref{DefI0}
with $q(\xi)=\p{\gamma\cdot\zeta}/\p{\rho\circ \psi}(\xi)$
is $L^2(\R^n;v)$--bounded.
Then we have
\[
\n{w(x)\zeta(D_x)\int^t_0
e^{i(t-\tau)a(D_x)}f(\tau,x)\,d\tau}_{L^2\p{\R_t\times\R^n_x}}
\leq C\n{v(x)f(t,x)}_{L^2\p{\R_t\times\R^n_x}}
\]
for all $f$ such that $\supp \FT_x f(t,\cdot)\subset\supp\gamma$,
where $a(\xi)=(\sigma\circ\psi)(\xi)$.
\end{theorem}
\medskip
\par
The proof of Theorem \ref{Th:reduction} is given in \cite{RS6},
and that of Theorem \ref{Th:redhom} is just a slight modification of it,
hence here we omit them.

As for the $L^2(\R^n;w)$--boundedness of the operator $I_{\psi,\gamma}$,
we have criteria for some special weight functions.
For $\ka\in\R$, let $L^2_\ka(\R^n)$, $\Dot{L}^2_\ka(\R^n)$ be
the set of measurable functions $f$ such that the norm
\[
\n{f}_{L^2_\ka(\R^n)}
=\p{\int_{\R^n}\abs{\langle x\rangle^\ka f(x)}^2\,dx}^{1/2},
\qquad
\n{f}_{\Dot{L}^2_\ka(\R^n)}
=\p{\int_{\R^n}\abs{|x|^\ka f(x)}^2\,dx}^{1/2}
\]
is finite, respectively.
Then we have the following mapping properties
(\cite[Theorems 4.2 and 4.3]{RS6}).
\par
\medskip
\begin{theorem}\label{Th:L2k}
Suppose $\ka\in\R$.
Assume that all the derivatives of entries of the 
$n\times n$ matrix
$\partial\psi$ and those of $\gamma$ are bounded.
Then the operators $I_{\psi,\gamma}$ and $I^{-1}_{\psi,\gamma}$ 
defined by
\eqref{DefI0} are $L^2_{\ka}(\R^n)$--bounded.
\end{theorem}
\begin{theorem}\label{Th:L'2k}
Let $\Gamma$, $\widetilde{\Gamma}\subset\R^n\setminus0$ be open cones.
Suppose $|\ka|< n/2$.
Assume $\psi(\lambda\xi)=\lambda\psi(\xi)$,
$\gamma(\lambda\xi)=\gamma(\xi)$ for all $\lambda>0$ and $\xi\in\Gamma$.
Then the operators $I_{\psi,\gamma}$ and $I^{-1}_{\psi,\gamma}$
defined by \eqref{DefI0} are $L^2_{\ka}(\R^n)$--bounded
and $\Dot{L}^2_{\ka}(\R^n)$--bounded.
\end{theorem}
\medskip
\par
We remark that the following result due to Kurtz and Wheeden 
\cite[Theorem 3]{KW} is essentially used to prove Theorem \ref{Th:L'2k}:
\medskip
\par
\begin{lemma}\label{Prop:wtbdd}
Suppose $|\ka|<n/2$.
Assume that $m(\xi)\in C^n(\R^n\setminus0)$
and all the derivative of $m(\xi)$ satisfies
$|\partial^\gamma m(\xi)|\leq C_\gamma
|\xi|^{-|\gamma|}$ for all $\xi\not=0$ and $|\gamma|\leq n$.
Then $m(D_x)$ is $L^2_{\ka}(\R^n)$ and 
$\Dot{L}^2_{\ka}(\R^n)$--bounded.
\end{lemma}
\medskip
\par

\section{Smoothing estimates for homogeneous dispersive equations}
\label{SECTION:main}

In author's paper \cite{RS6}, it is explained how to derive
smoothing estimates for general homogeneous dispersive equations
from model estimates
as an application of the canonical transformations described in Section
\ref{SECTION:canonical}.
We will repeat it here to help readers to understand the later part of
this paper concerning estimates for inhomogeneous equations.

Let us consider the solution
\[
u(t,x)=e^{ita(D_x)}\varphi(x)
\]
to the homogeneous equation
\[
\left\{
\begin{aligned}
\p{i\partial_t+a(D_x)}\,u(t,x)&=0\quad\text{in $\R_t\times\R^n_x$},\\
u(0,x)&=\varphi(x)\quad\text{in $\R^n_x$},
\end{aligned}
\right.
\]
where we always assume that 
function $a(\xi)$ is real-valued.
Let $a_m(\xi)\in C^\infty(\R^n\setminus0)$,
the {\it principal} part of $a(\xi)$, be 
a positively homogeneous function
of order $m$, that is, satisfy
$a_m(\lambda\xi)=\lambda^m a_m(\xi)$ for all $\lambda>0$ and $\xi\neq0$.
\par

First we consider the case that $a(\xi)$
has no lower order terms,
and assume that $a(\xi)$ is {\it dispersive}:
\medskip
\begin{equation}\tag{{\bf H}}
a(\xi)=a_m(\xi),\qquad\nabla a_m(\xi)\neq0 \quad(\xi\in\R^n\setminus0),
\end{equation}
\medskip
\par\noindent
where $\nabla=(\partial_1,\ldots,\partial_n)$ and
$\partial_j=\partial_{\xi_j}$.
A typical example is
$a(\xi)=a_m(\xi)=|\xi|^m$.
Especially, $a(\xi)=a_2(\xi)=|\xi|^2$
is the case of the Schr\"odinger equation.
\par
The following result (\cite[Theorem 5.1]{RS6})
is a generalisation of the one given by
Ben-Artzi and Klainerman \cite{BK} which treated the case
$a(\xi)=|\xi|^2$ and $n\geq 3$:
\par
\medskip
\begin{theorem}\label{M:H1}
Assume {\rm{(H)}}.
Suppose $n\geq 1$, $m>0$, and $s>1/2$.
Then we have
\begin{equation}\label{EQ:main1}
\n{\jp{x}^{-s}|D_x|^{(m-1)/2}e^{ita(D_x)}\varphi(x)}_{L^2\p{\R_t\times\R^n_x}}
\leq C\n{\varphi}_{L^2\p{\R^n_x}}.
\end{equation}
\end{theorem}
\medskip
\par
We review how to prove Theorem \ref{M:H1}.
The main idea is reducing them
to the special cases $a(D_n)=|D_n|^m,\,D_1|D_n|^{m-1}$,
where $D_x=(D_1,\ldots,D_n)$, by using Theorem \ref{Th:reduction}.
The following estimates (\cite[Theorem 3.1, Corollary 3.3]{RS6}) for them
act as model ones:
\par
\medskip
\par
\begin{proposition}\label{prop:basic}
Suppose $n=1$ and $m>0$.
Then we have
\[
\n{|D_x|^{(m-1)/2}e^{it|D_x|^m}\varphi(x)}_{L^2(\R_t)}
\leq C\n{\varphi}_{L^2(\R_x)}
\]
for all $x\in\R$.
Suppose $n=2$ and $m>0$.
Then we have
\[
\n{|D_y|^{(m-1)/2}e^{itD_x|D_y|^{m-1}}\varphi(x,y)}_{L^2(\R_t\times\R_y)}
\leq C\n{\varphi}_{L^2\p{\R^2_{x,y}}}
\]
for all $x\in\R$.
\end{proposition}
\begin{corollary}\label{Th:typeI}
Suppose $n\geq1$, $m>0$, and $s>1/2$.
Then we have
\begin{equation}\label{model:1}
\n{\jp{x_n}^{-s}|D_n|^{(m-1)/2}e^{it|D_n|^m}\varphi(x)}_{L^2(\R_t\times\R^n_x)}
\leq
 C\n{\varphi}_{L^2(\R_x^n)}.
\end{equation}
Suppose $n\geq2$, $m>0$, and $s>1/2$.
Then we have
\begin{equation}\label{model:2}
\n{\jp{x_1}^{-s}|D_n|^{(m-1)/2}e^{itD_1|D_n|^{m-1}}\varphi(x)}
_{L^2(\R_t\times\R^n_x)}
\leq
 C\n{\varphi}_{L^2(\R_x^n)}.
\end{equation}
\end{corollary}
\medskip
\par
We assume (H).
Let $\Gamma\subset\R^n\setminus0$ 
be a sufficiently small conic neighbourhood of
$e_n=(0,\ldots0,1)$, and take a cut-off function
$\gamma(\xi)\in C^\infty(\Gamma)$ which is 
positively homogeneous of
order $0$ and satisfies $\supp\gamma\cap \Sph^{n-1} 
\subset\Gamma\cap \Sph^{n-1}$.
By the microlocalisation and the rotation of the initial data
$\varphi$,
we may assume 
$\supp\widehat{\varphi}\subset\supp\gamma$.
The dispersive assumption $\nabla a_m(e_n)\neq0$ in this direction
implies the following two possibilities:
\medskip
\begin{description}
\item[(i)]
$\partial_n a_m(e_n)\neq0$.
Then, by Euler's identity $a_m(\xi)=(1/m)\nabla a_m(\xi)\cdot\xi$,
we have $a_m(e_n)\neq0$.
Hence, in this case, we may assume that
$a(\xi)(>0)$ and $\partial_n a(\xi)$ are bounded away from $0$
for $\xi\in\Gamma$.
\item[(ii)]
$\partial_n a_m(e_n)=0$.
Then there exits $j\neq n$ such that
$\partial_j a_m(e_n)\neq0$, say $\partial_1a_m(e_n)\neq0$.
Hence, in this case, we may assume
$\partial_1a(\xi)$ is bounded away from $0$
for $\xi\in\Gamma$.
We remark $a(e_n)=0$ by Euler's identity.
\end{description}
\medskip
\par
The estimate with the case $n=1$ is given by estimate \eqref{model:1}
in Corollary \ref{Th:typeI}.
In fact, we have $a(\xi)=a(1)|\xi|^m$ for $\xi>0$ in this case.
Hence we may assume $n\geq2$.
We remark that it is sufficient to show theorem with $1/2<s<n/2$
because the case $s\geq n/2$ is easily reduced to this case.
We will use the notation $\xi=(\xi_1,\ldots,\xi_n)$,
$\eta=(\eta_1,\ldots,\eta_n)$.
\par
In the case (i), we take
\begin{equation}\label{EQ:proofc1a1}
\sigma(\eta)=|\eta_n|^m,\quad
\psi(\xi)=(\xi_1,\ldots,\xi_{n-1},a(\xi)^{1/m}).
\end{equation}
Then we have $a(\xi)=\p{\sigma\circ\psi}(\xi)$ and
\begin{equation}\label{EQ:proofc1a2}
\det\partial\psi(\xi)
=
\begin{vmatrix}
E_{n-1}&0
\\
*&  (1/m)a(\xi)^{1/m-1} \partial_n a(\xi)
\end{vmatrix},
\end{equation}
where $E_{n-1}$ is the identity matrix of order $n-1$.
We remark that \eqref{infty} is satisfied since
$\det\partial\psi(e_n)=(1/m)a(e_n)^{1/m-1} \partial_n a(e_n)\neq0$.
By estimate \eqref{model:1} in Corollary \ref{Th:typeI},
we have
estimate \eqref{red} in Theorem \ref{Th:reduction} with $\sigma(D_x)=|D_n|^m$,
$w(x)=\jp{x}^{-s}$, and $\rho(\xi)=|\xi_n|^{(m-1)/2}$.
Note here the trivial inequality $\jp{x}^{-s}\leq\jp{x_n}^{-s}$.
If we take $\zeta(\xi)=|\xi|^{(m-1)/2}$,
then $q(\xi)=\gamma(\xi)\p{|\xi|/a(\xi)^{1/m}}^{(m-1)/2}$ defined by
\eqref{bdd} is a bounded function.
On the other hand, $I_{\psi,\gamma}$ is $L^2_{-s}$--bounded for $1/2<s<n/2$
by Theorem \ref{Th:L'2k}.
Hence, by Theorem \ref{Th:reduction}, we have estimate \eqref{org},
that is, estimate \eqref{EQ:main1}.
\par
In the case (ii), we take
\[
\sigma(\eta)=\eta_1|\eta_n|^{m-1},\quad
\psi(\xi)=\p{a(\xi)|\xi_n|^{1-m},\xi_2,\ldots,\xi_n}
\]
Then we have $a(\xi)=\p{\sigma\circ\psi}(\xi)$ and
\[
\det\partial\psi(\xi)
=
\begin{vmatrix}
\partial_1 a(\xi)|\xi_n|^{1-m}&*
\\
0&E_{n-1}
\end{vmatrix}.
\]
Since $\det\partial\psi(e_n)=\partial_1 a(e_n)\neq0$,
\eqref{infty} is satisfied.
Similarly to the case (i), the estimate for $\sigma(D_x)=D_1|D_n|^{m-1}$
is given by estimate \eqref{model:2} in Corollary \ref{Th:typeI},
which implies
estimate \eqref{EQ:main1} again by Theorem \ref{Th:reduction}.
\medskip
\par
As another advantage of the method explained here, we can also consider
the case that $a(\xi)$ has lower order terms,
and assume that $a(\xi)$ is dispersive in the following sense:
\medskip
\begin{equation}\tag{{\bf L}}
\begin{aligned}
&a(\xi)\in C^\infty(\R^n),\qquad \nabla a(\xi)\neq0 \quad(\xi\in\R^n),
\quad\nabla a_m(\xi)\neq0 \quad (\xi\in\R^n\setminus0), 
\\
&|\partial^\alpha\p{a(\xi)-a_m(\xi)}|\leq C_\alpha\abs{\xi}^{m-1-|\alpha|}
\quad\text{for all multi-indices $\alpha$ and all $|\xi|\geq1$}.
\end{aligned}
\end{equation}
or equivalently
\begin{equation}\tag{{\bf L}}
\begin{aligned}
&a(\xi)\in C^\infty(\R^n),\qquad 
|\nabla a(\xi)|\geq C\jp{\xi}^{m-1}\quad(\xi\in\R^n)\quad
\textrm{for some}\; C>0,
\\
&|\partial^\alpha\p{a(\xi)-a_m(\xi)}|\leq C_\alpha\abs{\xi}^{m-1-|\alpha|}
\quad\text{for all multi-indices $\alpha$ and all $|\xi|\geq1$}.
\end{aligned}
\end{equation}
The last lines of these assumptions simply amount to saying that
the principal part $a_m$ of $a$ is positively homogeneous
of order $m$ for $|\xi|\geq 1$.
\medskip
\par
The following result (\cite[Theorem 5.4]{RS6}) is also
derived from Corollary \ref{Th:typeI}:
\par
\medskip
\begin{theorem}\label{M:L4}
Assume {\rm (L)}.
Suppose $n\geq1$, $m>0$, and $s>1/2$.
Then we have
\begin{equation}\label{EQ:main4}
\n{\jp{x}^{-s}\jp{D_x}^{(m-1)/2}e^{ita(D_x)}\varphi(x)}
_{L^2\p{\R_t\times\R^n_x}}
\leq C\n{\varphi}_{L^2\p{\R^n_x}}.
\end{equation}
\end{theorem}
\medskip
\par
We review how to prove Theorem \ref{M:L4}.
We sometimes decompose the initial data $\varphi$ into the sum of
the {\it low frequency} part
$\varphi_{l}$ and the {\it high frequency} part $\varphi_{h}$, where
$\supp\widehat{\varphi_l}\subset\b{\xi:|\xi|< 2R}$
and $\supp\widehat{\varphi_h}\subset\b{\xi:|\xi|> R}$
with sufficiently large $R>0$.
Each part can be realised by multiplying $\chi(D_x)$ or $(1-\chi)(D_x)$
to $\varphi(x)$, hence to $u(t,x)$,
where $\chi\in C_0^\infty(\Rn)$ is an appropriate cut-off function.
For high frequency part,
the same argument as in the proof of Theorem \ref{M:H1} is valid.
(Furthermore, we can use Theorem \ref{Th:L2k} 
instead of Theorem \ref{Th:L'2k}
to assure the boundedness of $I_{\psi,\gamma}$, 
hence we need not assume
$n\geq2$.)
We show how to get the estimates for low frequency part.
Because of the compactness of it, we may assume
$\partial_ja(\xi)\neq0$ with some $j$, say $j=n$,
on a bounded set $\Gamma\subset\R^n$ and $\supp\widehat{\varphi}\subset\Gamma$.
Since we have $a(\xi)+c>0$ on $\Gamma$ with some constant $c>0$
and
\[
\n{\jp{x}^{-s}\jp{D_x}^{(m-1)/2}e^{ita(D_x)}\varphi}_{L^2(\R_t\times\R^n_x)}
=
\n{\jp{x}^{-s}\jp{D_x}^{(m-1)/2}
e^{it\p{a(D_x)+2c}}\varphi}_{L^2(\R_t\times\R^n_x)},
\]
we may assume $a(\xi)\geq c>0$ on $\Gamma$ without loss of generality.
We take a cut-off function
$\gamma(\xi)\in C_0^\infty(\Gamma)$,
and choose $\psi(\xi)$ and $\sigma(\eta)$ in the same way as
\eqref{EQ:proofc1a1}.
Assumption \eqref{infty} is also verified if we notice \eqref{EQ:proofc1a2}.
By estimate \eqref{model:1} in Corollary \ref{Th:typeI},
we have
estimate \eqref{red} in Theorem \ref{Th:reduction} with $\sigma(D_x)=|D_n|^m$,
$w(x)=\jp{x}^{-s}$ ($s>1/2$), and 
$\rho(\xi)=|\xi_n|^{(m-1)/2}$ as in
the proof of Theorem \ref{M:H1}.
If we take $\zeta(\xi)=\jp{\xi}^{(m-1)/2}$,
then $q(\xi)=\gamma(\xi)\p{\jp{\xi}/a(\xi)^{1/m}}^{(m-1)/2}$ 
defined by
\eqref{bdd} is a bounded function.
On the other hand, $I_{\psi,\gamma}$ is $L^2_{-s}$--bounded 
for all $s>1/2$
by Theorem \ref{Th:L2k}.
Hence, by Theorem \ref{Th:reduction}, we have estimate \eqref{org},
that is, estimate \eqref{EQ:main4}.
\par
\medskip
Finally, we introduce an intermediate assumption between (H) and (L),
and discuss what happens if we do not have
the condition $\nabla a(\xi)\neq0$:
\medskip
\begin{equation}\tag{{\bf HL}}
\begin{aligned}
&a(\xi)=a_m(\xi)+r(\xi),\quad
\nabla a_m(\xi)\neq0 \quad (\xi\in\R^n\setminus0),
\quad r(\xi)\in C^\infty(\R^n)
 \\
&|\partial^\alpha r(\xi)|\leq C\jp{\xi}^{m-1-|\alpha|}
\quad\text{for all multi-indices $\alpha$}.
\end{aligned}
\end{equation}
In view of the proof of Theorem \ref{M:L4}, we see that
Theorem \ref{M:H1} remains valid if we
replace assumption (H) by (HL) and functions $\varphi(x)$ in 
the estimates by its 
(sufficiently large) high frequency part $\varphi_h(x)$.
However we cannot control the low frequency part $\varphi_l(x)$,
and so have only the time local estimates on the whole.
We just put such a result (\cite[Theorem 5.6]{RS6})
below without its proof:
\par
\medskip
\begin{theorem}\label{Th:HL}
Assume {\rm (HL)}.
Suppose $n\geq1$, $m>0$, $s>1/2$, and $T>0$.
Then we have
\[
\int^T_0\n{\jp{x}^{-s}\jp{D_x}^{(m-1)/2}e^{ia(D_x)}}^2_{L^2(\R^n_x)}\,dt
\leq
 C\n{\varphi}_{L^2(\R^n)}^2,
\]
where $C>0$ is a constant depending on $T>0$.
\end{theorem}
\medskip
\par
We remark that Theorem \ref{M:L4} is the time global version
(that is, the estimate with $T=\infty$) of Theorem \ref{Th:HL},
and the extra assumption $\nabla a(\xi)\neq0$ is needed for that.
Since the assumption $\nabla a(\xi)\neq0$ for large $\xi$ 
is automatically satisfied by assumption (HL),
Theorem \ref{M:L4} means that the condition $\nabla a(\xi)\neq0$
for small $\xi$ assures the time global estimate.
In this sense, the low frequency part have a responsibility for
the time global smoothing. 
\bigskip
\par
\medskip

\section{Model estimates for inhomogeneous equations}
\label{SECTION:model-inh}
We now turn to deal with inhomogeneous equations, for which
we also have similar smoothing estimates. Such estimates are
necessary for nonlinear applications, and they can be
obtained by further developments of the presented methods.
Let us consider the solution
\[
u(t,x)
=-i\int^t_0e^{i(t-\tau)a(D_x)}f(\tau,x)\,d\tau
\]
to the equation
\[
\left\{
\begin{aligned}
\p{i\partial_t+a(D_x)}\,u(t,x)&=f(t,x)\quad\text{in $\R_t\times\R^n_x$},\\
u(0,x)&=0\quad\text{in $\R^n_x$}.
\end{aligned}
\right.
\]
We will give model estimates for it below,
where we write $x=(x_1,x_2,\ldots,x_n)\in\R^n$ and $D_x=(D_1,D_2\ldots,D_n)$.
We also write $x=x_1$, $D_x=D_1$ in the case $n=1$, 
and $(x,y)=(x_1,x_2)$, $(D_x,D_y)=(D_1,D_2)$ in the case $n=2$.
\par
\medskip
\begin{proposition}\label{Prop:inhom}
Suppose $n=1$ and $m>0$.
Let $a(\xi)\in C^\infty\p{\R\setminus0}$ be a real-valued
function which satisfies $a(\lambda\xi)=\lambda^m a(\xi)$ for
all $\lambda>0$ and $\xi\neq0$.
Then we have
\begin{equation}\label{eq:inhomdim1}
\n{a'(D_x)
\int^t_0e^{i(t-\tau)a(D_x)}f(\tau,x)\,d\tau}_{L^2(\R_t)}
\leq C\int_\R \n{f(t,x)}_{L^2(\R_t)}\,dx
\end{equation}
for all $x\in\R$.
Suppose $n=2$ and $m>0$.
Then we have
\begin{multline}\label{eq:inhomdim2}
\n{|D_x|^{m-1}\int^t_0e^{i(t-\tau)|D_x|^{m-1}D_y}
f(\tau,x,y)\,d\tau}_{L^2(\R_t\times\R_x)}
\\
\leq
C\int_\R \n{f(t,x,y)}_{L^2(\R_t\times\R_x)}\,dy
\end{multline}
for all $y\in\R$.
\end{proposition}
\begin{corollary}\label{Th:model-hom}
Suppose $n\geq1$, $m>0$, and $s>1/2$.
Let $a(\xi)\in C^\infty\p{\R\setminus0}$ be a real-valued
function which satisfies $a(\lambda\xi)=\lambda^m a(\xi)$ for
all $\lambda>0$ and $\xi\neq0$.
Then we have
\[
\n{\jp{x_n}^{-s}a'(D_n)
\int^t_0e^{i(t-\tau)a(D_n)}f(\tau,x)\,d\tau}
_{L^2(\R_t\times\R^n_x)}
\leq
 C\n{\jp{x_n}^{s}f(t,x)}_{L^2(\R_t\times\R_x^n)}.
\]
Suppose $n\geq2$, $m>0$, and $s>1/2$.
Then we have
\[
\n{\jp{x_1}^{-s}|D_n|^{m-1}\int^t_0e^{i(t-\tau)D_1|D_n|^{m-1}}f(\tau,x)\,d\tau}
_{L^2(\R_t\times\R^n_x)}
\leq
 C\n{\jp{x_1}^{s}f(t,x)}_{L^2(\R_t\times\R_x^n)}.
\]
\end{corollary}
\medskip
\par
Proposition \ref{Prop:inhom} with the case $n=1$
is a unification of the results
by Kenig, Ponce and Vega who treated the
cases $a(\xi)=\xi^2$ (\cite[p.258]{KPV3}), $a(\xi)=|\xi|\xi$
(\cite[p.160]{KPV4}),
and $a(\xi)=\xi^3$ (\cite[p.533]{KPV2}).
Corollary \ref{Th:model-hom} is a straightforward result of 
Proposition \ref{Prop:inhom} and Cauchy--Schwarz's inequality.
They act as model estimates for inhomogeneous equations
just like Proposition \ref{prop:basic} and Corollary \ref{Th:typeI}
do for homogeneous ones.
In \cite{RS6}, Corollary \ref{Th:typeI} is given straightforwardly
from the translation invariance of Lebesgue measure, by using a newly 
introduced method (comparison principle).

Since we unfortunately do not know the comparison principle
for inhomogeneous equations, we will give a direct proof to
Proposition \ref{Prop:inhom}.
Note that we have another expression 
of the solution to inhomogeneous equation
\[
\left\{
\begin{aligned}
\p{i\partial_t+a(D_x)}\,u(t,x)&=f(t,x)\quad\text{in $\R_t\times\R^n_x$},\\
u(0,x)&=0\quad\text{in $\R^n_x$},
\end{aligned}
\right.
\]
using the weak limit $R(\tau\pm i0)$ 
of the
resolvent $R(\tau\pm i\varepsilon)$ as $\varepsilon\searrow0$, where
$R(\lambda)=\p{a(D_x)-\lambda}^{-1}$:
\begin{equation}\label{resolv}
\begin{aligned}
u(t,x)
&=-i\int^t_0e^{i(t-\tau)a(D_x)}f(\tau,x)\,d\tau.
\\
&=\FT_\tau^{-1}R(\tau-i0)\FT_tf^+
+\FT_\tau^{-1}R(\tau+i0)\FT_tf^-
\end{aligned}
\end{equation}
(see Sugimoto \cite {Su1} and Chihara \cite{Ch}).
Here $\FT_t$ denotes the Fourier Transformation in $t$ and 
$\FT_\tau^{-1}$
its inverse, and $f^{\pm}(t,x)=f(t,x)Y(\pm t)$ is
the characteristic function $Y(t)$ of the set $\b{t\in\R:\,t>0}$.
\begin{proof}[Proof of Estimate \eqref{eq:inhomdim1}]
Let us use a variant of 
the argument of Chihara \cite[Section 4]{Ch}.
We set
$R(\lambda)=\p{a(D_x)-\lambda}^{-1}$
and show the estimate
\[
\abs{a^\prime (D_x)R(s\pm i0)g(x)}\leq C\int_\R|g(x)|\,dx,
\]
where $C>0$ is a constant independent of $s\in\R$, 
$x\in\R$ and $g\in L^1(\R)$.
Then, on account of the expression \eqref{resolv},
Plancherel's theorem, and Minkowski's inequality,
we have the desired result.
For this purpose, we consider the kernel
\[
k(s,x)=\FT^{-1}\br{a'(\xi)\p{a(\xi)-(s\pm i0)}^{-1}}(x)
\]
and show its uniform boundedness.
By the scaling argument, everything is reduced to show
the estimates
\[
\sup_{x\in \R}|k(\pm1,x)|\leq C
\qquad\text{and}\qquad \sup_{x\in \R}|k(0,x)|\leq C.
\]
By using an appropriate partition of
unity $\widehat{\phi}_1(\xi)+\widehat{\phi}_2(\xi)+\widehat{\phi}_3(\xi)=1$,
we split $k(\pm1, x)$ into the corresponding three parts
$k=k_1+k_2+k_3$,
where $\widehat{\phi}_1$ has its support near the origin,
$\widehat{\phi}_2$ near the point
$\xi^m=\pm1$, and $\widehat{\phi}_3$ away from these points.
The estimate for $k_1$ is trivial.
The other estimates are reduced to the boundedness of
\begin{equation}\label{heaviside}
k_0^\pm(x)=\FT^{-1}\br{\p{\xi\pm i0)}^{-1}}(x)=
\mp i\sqrt{2\pi}Y(\pm x).
\end{equation}
In fact, 
\[
k_2(\pm1,x)=\FT^{-1}\br{\p{\xi-(\alpha\pm i0)}^{-1}\widehat{\psi}(\xi)}(x)
=(e^{i\alpha x}k_0^\mp)*\psi(x)
\]
where $\alpha\in\R$ is a point which solves $a(\alpha)=\pm1$, and
\[
\widehat{\psi}(\xi)
=a'(\xi)\frac{\xi-\alpha}{a(\xi)-(\pm1)}
\widehat{\phi}_2(\xi)\in C^\infty_0(\R).\]
Furthermore, if we notice
\[
\frac{a'(\xi)}{a(\xi)-s}
=m\p{\frac{s}{(a(\xi)-s)\xi}+\frac1{\xi}},
\]
we have
\[
\frac1m
k_3(\pm1,x)=\pm\FT^{-1}\br{\frac{\widehat{\phi}_3(\xi)}
{(a(\xi)\mp1)\xi}}(x)
+k^\pm_0(x)-k^\pm_0*\p{\phi_1(x)+\phi_2(x)}.
\]
It is easy to deduce the estimates for $k_2$ and $k_3$.
It is also easy to verify
\[
\frac{a'(\xi)}{a(\xi)\pm i0}
=\frac{m}{\xi\pm i0}+c\delta
\]
with a constant $c$ and Dirac's delta function $\delta$,
and have the estimate for $k(0,x)$.
\end{proof}
\begin{proof}[Proof of Estimate \eqref{eq:inhomdim2}]
We set
$R(\lambda)=\p{|D_x|^{m-1}D_y-\lambda}^{-1}$
and show the estimate
\[
\n{|D_x|^{m-1}R(s\pm i0)g(x,y)}_{L^2(\R_x)}
\leq C\int\n{g(x,y)}_{L^2(\R_x)}\,dy,
\]
where $C>0$ is a constant independent of $s\in\R$,
$y\in\R$ and $g\in L^1(\R^2)$.
Then, by the expression \eqref{resolv}, Plancherel's theorem,
and Minkowski's inequality again,
we have the desired result.
\par
First we note, we may assume $\widehat{g}(\xi,\eta)=0$ for $\xi<0$.
Then we have
\begin{align*}
&|D_x|^{m-1}R(s\pm i0)g(x,y)
\\
=&\p{2\pi}^{-2}\int^\infty_0\int^\infty_{-\infty}
 e^{i(x\xi+y\eta)} |\xi|^{m-1}\p{|\xi|^{m-1}\eta-(s\pm i0)}^{-1}
 \widehat{g}(\xi,\eta)\,d\xi d\eta
\\
=&\p{2\pi}^{-2}\int^\infty_0\int^\infty_{-\infty}
 e^{ix\xi} |\xi|^{m-1}\p{|\xi|^{m-1}\eta-(s\pm i0)}^{-1}
 \widehat{g_y}(\xi,\eta)\,d\xi d\eta
\\
=&\p{2\pi}^{-2}\int^\infty_{-\infty}\int^\infty_{0}
 e^{ixb} \p{a-(s\pm i0)}^{-1}\widehat{g_y}(b,ab^{-(m-1)})\,dadb
\\
=&\p{2\pi}^{-2}\int^\infty_{-\infty}\int^\infty_{0}
 e^{ixb} \FT_a\br{\p{a-(s\pm i0)}^{-1}}
 \FT_a^{-1}\br{\widehat{g_y}(b,ab^{-(m-1)})}\,dadb
\\
=&\p{2\pi}^{-1}\int^\infty_{-\infty}\int^\infty_{0}
 e^{ixb} e^{-isa}k_0^\mp(-a)
 b^{m-1}\widetilde{g_y}(b,ab^{m-1})\,dadb,
\end{align*}
hence we have
\[
\FT_x\br{|D_x|^{m-1}R(s\pm i0)g(x,y)}(b)
=\int^\infty_{-\infty}
  e^{-isa}k_0^\mp(-a)
 b^{m-1}\widetilde{g_y}(b,ab^{m-1})\,da
\]
for $b\geq0$, and it vanishes for $b<0$.
Here $g_y(x,\,\cdot\,)=g(x,\,\cdot\,+y)$,
and $\widetilde{g_y}$ denotes its
partial Fourier transform with respect to the first variable.
We have also used here the change of variables
$a=\xi^{m-1}\eta$, $b=\xi$ and Parseval's formula.
Note that $\partial(a,b)/\partial(\xi,\eta)=b^{m-1}$ and
$k_0^\mp$ is a bounded function defined by \eqref{heaviside}.
Then we have the estimate
\begin{align*}
\abs{
\FT_x\br{|D_x|^{m-1}R(s\pm i0)g(x,y)}(b)
}
&\leq\sqrt{2\pi}
\int^\infty_{-\infty}
 \abs{b^{m-1}\widetilde{g_y}(b,ab^{m-1})}\,da
\\
&=\sqrt{2\pi}\int^\infty_{-\infty}
 \abs{\widetilde{g_y}(b,a)}\,da,
\end{align*}
and, by Plancherel's theorem and Minkowski's inequality, we have
\begin{align*}
\n{|D_x|^{m-1}R(s\pm i0)g(x,y)}_{L^2(\R_x)}
&\leq
\sqrt{2\pi}\int^\infty_{-\infty}
 \n{g_y(x,a)}_{L^2(\R_x)}\,da
\\
&=\sqrt{2\pi}\int^\infty_{-\infty}
 \n{g(x,y)}_{L^2(\R_x)}\,dy,
\end{align*}
which is the desired estimate.
\end{proof}
\medskip

\section{Smoothing estimates for dispersive inhomogeneous equations}
\label{SECTION:inhomogeneous}
Let us consider the inhomogeneous equation
\[
\left\{
\begin{aligned}
\p{i\partial_t+a(D_x)}\,u(t,x)&=f(t,x)\quad\text{in $\R_t\times\R^n_x$},\\
u(0,x)&=0\quad\text{in $\R^n_x$},
\end{aligned}
\right.
\]
where we always assume that 
function $a(\xi)$ is real-valued.
Let the principal part $a_m(\xi)\in C^\infty(\R^n\setminus0)$,
be a positively homogeneous function
of order $m$.
Recall the dispersive conditions we used in Section \ref{SECTION:main}:
\medskip
\begin{equation}\tag{{\bf H}}
a(\xi)=a_m(\xi),\qquad\nabla a_m(\xi)\neq0 \quad(\xi\in\R^n\setminus0),
\end{equation}
\begin{equation}\tag{{\bf L}}
\begin{aligned}
&a(\xi)\in C^\infty(\R^n),\qquad \nabla a(\xi)\neq0 \quad(\xi\in\R^n),
\quad\nabla a_m(\xi)\neq0 \quad (\xi\in\R^n\setminus0), 
\\
&|\partial^\alpha\p{a(\xi)-a_m(\xi)}|\leq C_\alpha\abs{\xi}^{m-1-|\alpha|}
\quad\text{for all multi-indices $\alpha$ and all $|\xi|\geq1$}.
\end{aligned}
\end{equation}
\medskip
The following is a counterpart of Theorem \ref{M:H1}
which treated homogeneous equations:
\par
\medskip
\begin{theorem}\label{Th:inhom}
Assume {\rm (H)}.
Suppose $m>0$ and $s>1/2$.
Then we have
\begin{equation}\label{eq:inhom:1}
\n{\jp{x}^{-s}|D_x|^{m-1}\int^t_0e^{i(t-\tau)a(D_x)}f(\tau,x)\,d\tau}_
{L^2(\R_t\times\R^n_x)}
\leq
 C\n{\jp{x}^s f(t,x)}_{L^2(\R_t\times\R^n_x)}
\end{equation}
in the case $n\geq2$, and
\begin{equation}\label{eq:inhom:2}
\n{\jp{x}^{-s}a'(D_x)\int^t_0e^{i(t-\tau)a(D_x)}f(\tau,x)\,d\tau}_
{L^2(\R_t\times\R_x)}
\leq
 C\n{\jp{x}^s f(t,x)}_{L^2(\R_t\times\R_x)}
\end{equation}
in the case $n=1$.
\end{theorem}
\medskip
\par
Chihara \cite{Ch} proved Theorem \ref{Th:inhom} with $m>1$
under the assumption \rm{(H)}.
As was pointed out in \cite[p.1958]{Ch},
we cannot replace $a'(D_x)$ by $|D_x|^{m-1}$ in
estimate \eqref{eq:inhom:2} for the case
$n=1$, but there is another explanation for this obstacle.
If we decompose
$f(t,x)=\chi_+(D_x)f(t,x)+\chi_-(D_x)f(t,x)$,
where $\chi_{\pm}(\xi)$ is a characteristic function of the 
set $\b{\xi\in\R\,:\,\pm\xi\geq0}$,
then we easily obtain
\begin{align*}
&\n{\jp{x}^{-s}|D_x|^{(m-1)/2}\int^t_0e^{i(t-\tau)a(D_x)}f(\tau,x)\,d\tau}_
{L^2(\R_t\times\R_x)}
\\
\leq
 &C\Bigl(
\n{\jp{x}^s |D_x|^{-(m-1)/2}f_+(t,x)}_{L^2(\R_t\times\R_x)}
 +\n{\jp{x}^s |D_x|^{-(m-1)/2}f_-(t,x)}_{L^2(\R_t\times\R_x)}
\Bigr)
\end{align*}
from Theorem \ref{Th:inhom}.
But we cannot justify the estimate
\[
 \n{\jp{x}^s|D_x|^{-(m-1)/2} f_\pm(t,x)}_{L^2(\R_t\times\R_x)}
\leq
 C \n{\jp{x}^s|D_x|^{-(m-1)/2} f(t,x)}_{L^2(\R_t\times\R_x)}
\]
for $s>1/2$ by Lemma \ref{Prop:wtbdd}
because it requires $s<n/2$ and it is impossible for $n=1$.
\par
As a counterpart of Theorem \ref{M:L4}, we have
\par
\medskip
\begin{theorem}\label{Th:inhom2}
Assume {\rm (L)}.
Suppose $n\geq1$, $m>0$, and $s>1/2$.
Then we have
\begin{equation}\label{eq:inhom:3}
\n{\jp{x}^{-s}\jp{D_x}^{m-1}\int^t_0e^{i(t-\tau)a(D_x)}f(\tau,x)\,d\tau}_
{L^2(\R_t\times\R^n_x)}
\leq
 C\n{\jp{x}^s f(t,x)}_{L^2(\R_t\times\R^n_x)}.
\end{equation}
\end{theorem}
\medskip
\par
The following result is a straightforward consequence
of Theorem \ref{Th:inhom2} and the $L^2_s$--boundedness
of $|D_x|^{(m-1)/2}\jp{D_x}^{-(m-1)/2}$
with $(1/2<)s<n/2$ and $m\geq1$ (which is assured by Lemma \ref{Prop:wtbdd}):
\par
\medskip
\begin{corollary}\label{Cor:inhom2}
Assume {\rm (L)}.
Suppose $n\geq2$, $m\geq1$, and $s>1/2$.
Then we have
\[
\n{\jp{x}^{-s}|D_x|^{m-1}\int^t_0e^{i(t-\tau)a(D_x)}f(\tau,x)\,d\tau}_
{L^2(\R_t\times\R^n_x)}
\leq
 C\n{\jp{x}^s f(t,x)}_{L^2(\R_t\times\R^n_x)}.
\]
\end{corollary}
\medskip
\par
We remark that the same argument of canonical transformations
as used for homogeneous equations in
Section \ref{SECTION:main} works for inhomogeneous ones, as well.
That is, the proofs of Theorems \ref{Th:inhom} and \ref{Th:inhom2}
are carried out by reducing them to model estimates in
Corollary \ref{Th:model-hom}.
We omit the details because the argument is essentially the same,
but we just remark that we use Theorem \ref{Th:redhom} instead of
Theorem \ref{Th:reduction}.

The following is a counterpart of Theorem \ref{Th:HL}:
\par
\medskip
\begin{theorem}\label{Th:HLinh}
Assume {\rm (HL)}.
Suppose $n\geq1$, $m>0$, $s>1/2$, and $T>0$.
Then we have
\begin{multline*}
\int^T_0\n{\jp{x}^{-s}\jp{D_x}^{m-1}
\int^t_0e^{i(t-\tau)a(D_x)}f(\tau,x)\,d\tau}^2_{L^2(\R^n_x)}\,dt
\\
\leq
 C\int^T_0\n{\jp{x}^sf(t,x)}^2_{L^2(\R^n_x)}\,dt,
\end{multline*}
where $C>0$ is a constant depending on $T>0$.
\end{theorem}
\medskip
\begin{proof}
By multiplying $\chi(D_x)$ and $(1-\chi)(D_x)$ to $f(t,x)$,
we decompose it into the sum 
of low frequency part and high frequency part,
where $\chi(\xi)$ is an appropriate cut-off function.
As in the proof of Theorem \ref{M:L4},
the estimate for the high frequency part can be reduced to
Corollary \ref{Th:model-hom}
by using Theorem \ref{Th:redhom} instead of Theorem \ref{Th:reduction},
together with the boundedness result Theorem \ref{Th:L2k}.
Here we note that, for $t\in[0,T]$,
\[
\int^t_0 e^{i(t-\tau)a(D_x)}f(\tau,x)\,d\tau
=
\int^t_0 e^{i(t-\tau)a(D_x)}\chi_{[0,T]}(\tau)f(\tau,x)\,d\tau,
\]
where $\chi_{[0,T]}$ denotes the characteristic 
function of the interval $[0,T]$.
The estimate for the low frequency part is trivial.
In fact, if $\supp_\xi \FT_xf(t,\xi)\subset
\br{\xi;|\xi|\leq R}$, we have
\begin{align*}
&\int^T_0
\n{\jp{x}^{-s}\jp{D_x}^{m-1}\int^t_0e^{i(t-\tau)a(D_x)}
 f(\tau,x)\,d\tau}^2_{L^2(\R^n_x)}
\,dt
\\
\leq
&\int^T_0
\p{
\int^T_0
\n{
\jp{D_x}^{m-1}e^{i(t-\tau)a(D_x)}f(\tau,x)
}_{L^2(\R^n_x)} \,d\tau
}^2
\,dt
\\
\leq &CT^2\jp{R}^{2(m-1)}
 \int^T_0\n{\jp{x}^sf(t,x)}^2_{L^2(\R^n_x)}\,dt.
\end{align*}
by Plancherel's theorem.
\end{proof}
\par
\medskip
If we combine Theorem \ref{Th:inhom} with Theorem \ref{M:H1},
we have a result for the equation
\begin{equation}\label{eq:inhom2}
\left\{
\begin{aligned}
\p{i\partial_t+a(D_x)}\,u(t,x)&=f(t,x)\quad\text{in $\R_t\times\R^n_x$},\\
u(0,x)&=\varphi(x)\quad\text{in $\R^n_x$}.
\end{aligned}
\right.
\end{equation}
\par
\medskip
\begin{corollary}\label{Cor:inhom+hom}
Assume {\rm (H)}.
Suppose $m>0$ and $s>1/2$.
Then the solution $u$ to equation \eqref{eq:inhom2}
satisfies
\begin{multline*}
\n{\jp{x}^{-s}|D_x|^{-(m-1)/2} a'(D_x)u(t,x)}
_{L^2(\R_t\times\R_x)}
\\
\leq
 C\p{\n{\varphi}_{L^2(\R)}
+\n{\jp{x}^s|D_x|^{-(m-1)/2}f(t,x)}_{L^2(\R_t\times\R_x)}}
\end{multline*}
in the case $n=1$, and
\begin{multline*}
\n{\jp{x}^{-s}|D_x|^{(m-1)/2}u(t,x)}
_{L^2(\R_t\times\R^n_x)}
\\
\leq
 C\p{\n{\varphi}_{L^2(\R^n)}
+\n{\jp{x}^s|D_x|^{-(m-1)/2}f(t,x)}_{L^2(\R_t\times\R^n_x)}}
\end{multline*}
in the case $n\geq2$.
\end{corollary}
\medskip
\par
If we combine Theorem \ref{Th:inhom2} with Theorem \ref{M:L4},
we have the following:
\par
\medskip
\begin{corollary}\label{Cor:inhom+hom2}
\medskip
Assume {\rm (L)}.
Suppose $n\geq1$, $m>0$, and $s>1/2$.
Then the solution $u$ to equation \eqref{eq:inhom2}
satisfies
\begin{multline*}
\n{\jp{x}^{-s}\jp{D_x}^{(m-1)/2}u(t,x)}
_{L^2(\R_t\times\R_x)}
\\
\leq
 C\p{\n{\varphi}_{L^2(\R)}
+\n{\jp{x}^s\jp{D_x}^{-(m-1)/2}f(t,x)}_{L^2(\R_t\times\R_x)}}.
\end{multline*}
\end{corollary}
\medskip
\par
If we combine Theorem \ref{Th:HLinh} with Theorem \ref{Th:HL},
we have the following:
\par
\medskip
\begin{corollary}\label{Cor:HL}
Assume {\rm (HL)}.
Suppose $n\geq1$, $m>0$, $s>1/2$, and $T>0$.
Then the solution $u$ to equation \eqref{eq:inhom2} satisfies
\begin{multline*}
\int^T_0\n{\jp{x}^{-s}\jp{D_x}^{(m-1)/2}u(t,x)}^2_{L^2(\R^n_x)}\,dt
\\
\leq
 C\p{\n{\varphi}_{L^2(\R^n)}^2
+\int^T_0\n{\jp{x}^s\jp{D_x}^{-(m-1)/2}f(t,x)}^2_{L^2(\R^n_x)}\,dt},
\end{multline*}
where $C>0$ is a constant depending on $T>0$.
\end{corollary}
\medskip
\par
Corollary \ref{Cor:HL} is an extension of the result
by Hoshiro \cite{Ho2}, which treated the case that $a(\xi)$ 
is a polynomial.
The proof relied on Mourre's method, which is known in spectral
and scattering theories. Here we use the argument of canonical 
transformations, extending the result and simplifying the proof.



\begin{thebibliography}{1}



%


\bibitem[BD1]{BD1}
M.~Ben-Artzi and A.~Devinatz, 
{\it The limiting absorption principle for partial 
differential operators}, 
Mem. Amer. Math. Soc. {\bf 66} (1987). 

\bibitem[BD2]{BD2} M.~Ben-Artzi and A.~Devinatz,
 {\it Local smoothing and convergence properties of Schr\"{o}dinger
 type equations},
 J. Funct. Anal. {\bf 101} (1991), 231--254.

\bibitem[BK]{BK} M.~Ben-Artzi and S.~Klainerman,
 {\it Decay and regularity for the Schr\"{o}dinger equation},
 J. Analyse Math. {\bf 58} (1992),
 25--37.
 
 

%
%

%
%

\bibitem[Ch]{Ch} H.~Chihara,
  {\it Smoothing effects of dispersive pseudodifferential equations},
  Comm. Partial Differential Equations
  {\bf 27} (2002), 1953--2005.







\bibitem[CS]{CS} P.~Constantin and J.~C.~Saut,
 {\it Local smoothing properties of dispersive equations},
 J. Amer. Math. Soc. {\bf 1} (1988),
 413--439.

%
%
%
%
%





%


\bibitem[Ho1]{Ho1} T.~Hoshiro, {\it Mourre's method and smoothing
 properties of dispersive equations}, Comm. Math. Phys. 
 {\bf 202} (1999), 255--265. 

\bibitem[Ho2]{Ho2} T.~Hoshiro,
 {\it Decay and regularity for dispersive equations with
     constant coefficients},
 J. Anal. Math. {\bf 91} (2003), 211--230.

%


\bibitem[Ka2]{Ka2} T.~Kato, 
{\it On the Cauchy problem for the (generalized) 
Korteweg-de Vries equation}, 
Studies in applied mathematics, 93--128, 
Adv. Math. Suppl. Stud., 8, Academic Press, New York, 1983.

\bibitem[KY]{KY} T.~Kato and K.~Yajima,
 {\it Some examples of smooth operators and the 
 associated smoothing effect},
 Rev. Math. Phys. {\bf 1} (1989),
 481--496.

\bibitem[KPV1]{KPV1} C.~E.~Kenig, G.~Ponce and L.~Vega,
{\it Oscillatory integrals and regularity of dispersive equations},
Indiana Univ. Math. J. {\bf 40} (1991), 33--69.

\bibitem[KPV2]{KPV2} C.~E.~Kenig, G.~Ponce and L.~Vega,
{\it Well-posedness and scattering results for the generalized
 Korteweg-de Vries equation via the contraction principle},
Comm. Pure Appl. Math. {\bf 46} (1993), 527--620.

\bibitem[KPV3]{KPV3} C.~E.~Kenig, G.~Ponce and L.~Vega,
{\it Small solutions to nonlinear Schr\"odinger equations},
Ann. Inst. H. Poincar\'e Anal. Non Lin\'eaire {\bf 10} 
(1993), 255--288.

\bibitem[KPV4]{KPV4} C.~E.~Kenig, G.~Ponce and L.~Vega,
{\it On the generalized Benjamin-Ono equation},
Trans. Amer. Math. Soc. {\bf 342} (1994), 155--172.

\bibitem[KPV5]{KPV5} C.~E.~Kenig, G.~Ponce and L.~Vega, 
{\it On the Zakharov and Zakharov-Schulman systems}, 
J. Funct. Anal. {\bf 127} (1995), 204--234. 


\bibitem[KoSa]{KoSa} H.~Koch and J.~C.~Saut,
{\it Local smoothing and local solvability for third
order dispersive equations},
preprint.


\bibitem[KW]{KW} D.~S.~Kurtz and R.~L.~Wheeden,
 {\it Results on weighted norm inequalities for multipliers},
 Trans. Amer. Math. Soc. {\bf 255} (1979), 343--362.

\bibitem[La]{La} C.~Laurey, 
{\it The Cauchy problem for a third order nonlinear 
Schr\"odinger equation},
Nonlinear Anal. {\bf 29} (1997), 121--158. 

\bibitem[LP]{LP}  F.~Linares and G.~Ponce,
 {\it On the Davey-Stewartson systems},
 Ann. Inst. H. Poincar\'e Anal. Non Lin\'eaire {\bf 10} (1993), 523--548.

%


 
%
%
%

\bibitem[RS1]{RS1}
M.~Ruzhansky and M.~Sugimoto,
{\it A new proof of global smoothing 
estimates for dispersive equations},
Advances in pseudo-differential operators, 65--75,
Oper. Theory Adv. Appl., {\bf 155}, Birkh\"auser, Basel, 2004.



\bibitem[RS2]{RS4} M.~Ruzhansky and M.~Sugimoto,
 {\it Global smoothing estimates for dispersive equations
 with non-polynomial symbols},
 Proceedings of The 12th International Conference on
    Finite or Infinite Dimensional Complex Analysis and Applications,
    Kyushu University Press, Fukuoka.



\bibitem[RS3]{RS6}   
 M.~Ruzhansky and M.~Sugimoto, 
 {\it Smoothing properties of evolution equations via canonical
 transforms and comparison principle}, Proc. London Math. Soc. 
 {\bf 105} (2012), 393--423. 
   
 \bibitem[RS4]{RScmp}   
 M.~Ruzhansky and M.~Sugimoto, 
 {\it Structural resolvent estimates and derivative nonlinear
 Schr\"odinger equations}, Comm. Math. Phys.
{\bf 314} (2012), 281--304.


%

\bibitem[Si]{Si} B.~Simon, {\it Best constants in some operator
  smoothness estimates}, J. Funct. Anal. {\bf 107} (1992), 66--71.

\bibitem[Sj]{Sj} P.~Sj\"olin,
 {\it Regularity of solutions to the Schr\"odinger equation},
 Duke Math. J. {\bf 55} (1987), 699--715.






\bibitem[Su1]{Su1} M.~Sugimoto,
 {\it Global smoothing properties of generalized 
 Schr\"odinger equations},
 J. Anal. Math. {\bf 76} (1998), 191--204.

\bibitem[Su2]{Su2} M.~Sugimoto,
 {\it A Smoothing property of Schr\"odinger equations 
 along the sphere},
 J. Anal. Math. {\bf 89} (2003), 15--30.

%

\bibitem[V]{V} L.~Vega,
{\it Schr\"odinger equations: Pointwise convergence to the
initial data}, Proc. Amer. Math. Soc. {\bf 102} (1988),
874--878.


\bibitem[Wa1]{Wa1} B.~G.~Walther,
 {\it A sharp weighted $L^2$-estimate for the solution to the
 time-dependent Schr\"odinger equation},
 Ark. Mat. {\bf 37} (1999), 381--393.

\bibitem[Wa2]{Wa2} B.~G.~Walther,
 {\it Regularity, decay, and best constants for dispersive equations},
 J. Funct. Anal. {\bf 189} (2002), 325--335.

%

\end{thebibliography}
\end{document}